\documentclass{article}
\usepackage{amsmath}
\usepackage{amsfonts,amsthm}
\usepackage{amssymb,enumerate,enumitem,verbatim}
\usepackage[pdftex]{graphicx}

\newcounter{myfootnote}[page]

\newtheorem{lemma}{Lemma}[section]
\newtheorem{corollary}[lemma]{Corollary}
\newtheorem{theorem}[lemma]{Theorem}
\newtheorem{prop}[lemma]{Proposition}

\theoremstyle{definition}
\newtheorem*{acknowledgements}{Acknowledgements}
\newtheorem*{rmk}{Remark}
\newtheorem*{defn}{Definition}
\newtheorem*{conjecture}{Conjecture}

\newtheorem*{claim}{Claim}

\theoremstyle{remark}

\newtheorem*{example}{Example}
\newtheorem*{nte}{Note}

\global\long\def\a{\alpha}
\global\long\def\b{\beta}

\global\long\def\d{\delta}

\global\long\def\e{\varepsilon}

\global\long\def\N{\mathbb{N}}

\global\long\def\GG{\mathcal{G}}

\global\long\def\AA{\mathcal{A}}
\global\long\def\BB{\mathcal{B}}
\global\long\def\CC{\mathcal{C}}
\global\long\def\DD{\mathcal{D}}

\global\long\def\UU{\mathcal{U}}
\global\long\def\VV{\mathcal{V}}

\global\long\def\TT{\mathcal{T}}

\global\long\def\re{\begin{rmk}}
\global\long\def\mark{\end{rmk}}
\global\long\def\ex{\begin{example}}
\global\long\def\ple{\end{example}}
\global\long\def\no{\begin{nte}}
\global\long\def\ted{\end{nte}}
\global\long\def\en{\begin{compactenum}}
\global\long\def\um{\end{compactenum}}
\global\long\def\li{\begin{compactitem}}
\global\long\def\st{\end{compactitem}}
\global\long\def\de{\begin{defn}}
\global\long\def\fn{\end{defn}}
\global\long\def\cor{\begin{corollary}}
\global\long\def\ary{\end{corollary}}
\global\long\def\lem{\begin{lemma}}
\global\long\def\ma{\end{lemma}}
\global\long\def\arr{\begin{array}}
\global\long\def\ay{\end{array}}
\global\long\def\pr{\begin{proof}}
\global\long\def\oof{\end{proof}}

\newif\ifdraft
\draftfalse % or \drafttrue\draftfalse

\newif\ifdone
\donefalse
\title{Sharp threshold for embedding combs and other spanning trees in random graphs}

\author{Richard Montgomery\footnote{Department of Pure Mathematics and Mathematical Statistics, Centre for Mathematical Sciences, Wilberforce Road, Cambridge, CB3 0WB, UK. r.h.montgomery@dpmms.cam.ac.uk}}

\begin{document}

\maketitle
\begin{abstract}
When $k|n$, the tree $\mathrm{Comb}_{n,k}$ consists of a path containing $n/k$ vertices, each of whose vertices has a disjoint path length $k-1$ beginning at it. We show that, for any $k=k(n)$ and $\e>0$, the binomial random graph $\GG(n,(1+\e)\log n/ n)$ almost surely contains $\mathrm{Comb}_{n,k}$ as a subgraph. This improves a recent result of Kahn, Lubetzky and Wormald. We prove a similar statement for a more general class of trees containing both these combs and all bounded degree spanning trees which have at least $\e n/ \log^9n$ disjoint bare paths length $\left\lceil\log^9 n\right\rceil$.

We also give an efficient method for finding large expander subgraphs in a binomial random graph. This allows us to improve a result on almost spanning trees by Balogh, Csaba, Pei and Samotij.
\end{abstract}

\section{Introduction}\label{intro}

Given a tree $T$ with $n$ vertices and maximum degree at most $\Delta$, for what range of $p$ are we likely to find a copy of $T$ in the binomial random graph $\GG(n,p)$? Around twenty years ago, Kahn proposed the following natural conjecture~\cite{KLW14b}.

\begin{conjecture}
For every fixed $\Delta>0$, there is some constant $C$ such that if $T$ is a tree on $n$ vertices with maximum degree $\Delta$, then the random graph $\GG(n,C\log n/ n)$ almost surely contains a copy of $T$.
\end{conjecture}
\noindent As we expect isolated vertices when $p<\log n/ n$, this conjecture would be tight up to the constant.

%\footnote{More formally, given any sequence of such trees with $|T_n|=n$, $\P(T_n\subset\GG(n,C\log n/ n)\to 1$.}

Alon, Krivelevich and Sudakov~\cite{AKS07} showed that, for every $\e>0$ and $\Delta\in\N$, there is some $c=c(\e,\Delta)$ for which the random graph $\GG(n,c/n)$ almost surely contains a copy of every tree on at most $(1-\e)n$ vertices with maximum degree at most $\Delta$. Taking $\TT(n,\Delta)$ to be the class of all trees on $n$ vertices with maximum degree at most $\Delta$, we say such a graph is $\TT((1-\e)n,\Delta)$-universal. Alon, Krivelevich and Sudakov used their result for these \emph{almost spanning trees} to demonstrate that the above conjecture for spanning trees is true for trees which have at least $\a n$ leaves, for any fixed $\a>0$.

Krivelevich~\cite{MK10} showed that, for every $\e>0$, if $T\in \TT(n,\Delta)$, then $\GG(n,n^{-1+\e})$ almost surely contains a copy of $T$, as well as considering the same question with larger, non-constant, $\Delta$. The author~\cite{selflove1} recently proved that, if $T\in\TT(n,\Delta)$, then $\GG(n,\Delta\log^5 n/ n)$ almost surely contains a copy of $T$.

Hefetz, Krivelevich and Szab\'{o}~\cite{HKS12} proved that, for each $\e,\a>0$, if $T\in \TT(n,\Delta)$ has at least $\a n$ leaves, then $\GG(n,(1+\e)\log n/n)$ almost surely contains a copy of $T$. As noted above, the threshold for $G(n,p)$ to be connected demonstrates that the probability they used is tight up to a factor of $(1+\e)$. Hefetz, Krivelevich and Szab\'{o} were also able to prove that such a random graph almost surely contains all the trees in $\TT(n,\Delta)$ which contain a path of length at least $\a n$ whose interior vertices all have degree 2 in $T$. Such a path is known as a \emph{bare path}.

As early as the statement of the conjecture above, the study of the following specific trees, known as \emph{combs}, was suggested~\cite{KLW14b}. When $k|n$, the tree $\mathrm{Comb}_{n,k}$ consists of a path $P$ of length $n/k-1$, and $n/k$ disjoint paths of length $k-1$, each starting at a vertex of $P$. The conjecture was recently proved for such trees by Kahn, Lubetzky and Wormald~\cite{KLW14a,KLW14b}. Here, we show that $p=(1+\e)\log n/ n$ is sufficient to almost surely find a copy of a tree $T$ in $\GG(n,p)$ if $T$ belongs to a wider class of trees which contains the combs. 

We say a path $P$ in a tree $T$ is a \emph{tooth} if it is a bare path, one of whose end vertices is a leaf. The tree $\mathrm{Comb}_{n,k}$ has $n/k$ teeth length $k$, and thus is covered by the following theorem.

\begin{theorem} \label{combshort}
Let $\a,\e>0$ and $\Delta\in\N$ be fixed. Let $k\geq 10$ and suppose $T\in \TT(n,\Delta)$ has at least $\a n/k$ teeth length $k$. Then, almost surely, the random graph $\GG(n,(1+\e)\log n / n)$ contains a copy of $T$.
\end{theorem}

When $k$ is fixed, the trees in Theorem \ref{combshort} have at least $\a n/k$ leaves and are covered by the work of Hefetz, Krivelevich and Szab\'o. The methods used by Kahn, Lubetzky and Wormald could be used to embed the trees considered in Theorem \ref{combshort}, but would require the probability to increase with $\a$. By combining the methods used to prove Theorem \ref{combshort} with methods in~\cite{selflove1}, we extend the class of trees $T$ for which we can almost surely find a copy of $T$ in $\GG(n,(1+\e)\log n/ n)$.

\begin{theorem} \label{comblong}
Let $\a,\e,\Delta>0$. The random graph $\GG(n,(1+\e)\log n / n)$ almost surely contains a copy of every tree $T\in \TT(n,\Delta)$ which has at least $\a n/\log^9 n$ disjoint bare paths with length $\left\lceil \log^9 n \right\rceil$.
\end{theorem}

In fact, the proof for Theorem \ref{comblong} requires of a graph $T$ only that it has bounded degree and is acyclic once the bare paths of length $\left\lceil \log^9 n\right\rceil$ are removed. Thus, the following holds.

\begin{theorem}\label{cycles}
For each $\e>0$, if $k|n$ and $k\geq 2\log^9 n$, then almost surely $\GG(n,(1+\e)\log n/n)$ can be covered by $n/k$ cycles of length $k$.
\end{theorem}

In the range of $k$ specified this improves a result of Lubetzky, Kahn and Wormald factoring the binomial random graph into cycles~\cite{KLW14a}. As such a factorisation into cycles requires $\GG(n,p)$ to have no isolated vertices, this result is tight. 

Lubetzky, Kahn and Wormald's approach to embedding $\mathrm{Comb}_{n,k}$ in the region $k\geq \kappa\log n$, for some fixed $\kappa$, used their almost sure factorisation of $\GG(n,p)$ into cycles length $k$. They show this factorisation exists using a delicate probabilistic argument. We take a very different, constructive, approach, that only uses simple probability to derive expansion properties of a random graph, which are then used to find the required spanning trees.

As part of our efforts to give a clear presentation of the proofs of Theorems \ref{combshort} and \ref{comblong}, we make use of a simple method for finding, in graphs with expansion properties for large sets, large subgraphs with expansion properties for small sets. As an example of this, we make a small improvement on the probability required to embed almost spanning trees. Balogh, Csaba, Pei and Samotij~\cite{BCPS10} used a theorem of Haxell~\cite{PH01} to show that the random graph $\GG(n,c/n)$ is almost surely $\TT((1-\e)n,\Delta)$-universal if $c\geq \max\{1000\Delta\log(20 \Delta), 30(\Delta/\e)\log (4e/\e)\}$. This improved the value of $c$ used by Alon, Krivelevich and Sudakov for the same result. We simplify the proof in~\cite{BCPS10}, giving a small improvement to the value of $c$ required.
\begin{theorem}\label{almostprob} Let $d\geq 2$ and $0<\e<1/2$. If $c\geq (30\Delta/\e)\log(4e/\e)$, then the random graph $\GG(n,c/n)$ is almost surely $\TT((1-\e)n,\Delta)$-universal.
\end{theorem}
As noted by Alon, Krivelevich and Sudakov~\cite{AKS07}, the constant $c$ cannot be reduced beneath $c_0\Delta\log(1/\e)$ for some absolute constant $c_0$. Hence, the dependence of $c$ on $\Delta$ in Theorem \ref{almostprob} is correct.

In Section \ref{tools} we will cover some simple probabilistic results, as well as collect tools from other work. In Section \ref{almostsec} we will prove Theorem \ref{almostprob}. In Section \ref{rotate} we introduce $(l,\gamma)$-connectors and sketch an embedding of the comb in the random graph $\GG(n,\log^2 n/n)$. Sections \ref{secshort} and \ref{seclong} prove Theorems \ref{combshort} and \ref{comblong} respectively.

\section{Preliminaries}\label{tools}

\subsection{Notation}
For a graph $G$, $V(G)$ will be the vertex set of $G$ and $|G|=|V(G)|$. Where $W\subset V(G)$, $G[W]$ is the subgraph of $G$ induced on the vertices of $W$. The set of neighbours of a vertex $v$ is denoted by $N(v)$, and the neighbourhood of a vertex set $A\subset V(G)$ by $N(A)=(\cup_{v\in A}N(v))\setminus A$. Where multiple graphs are used, we use $N_G(v)$ for the neighbourhood of a vertex $v$ in the graph $G$. We use $N(A,B)$ to refer to the set of neighbours of $A$ in $B$, that is $N(A)\cap B$. We take $d_G(x, A)=|N_G(x,A)|$, and let
\[
d_G(A,B)=\sum_{x\in A}d_G(x,B).
\]

We say a path with $l$ vertices has \emph{length} $l-1$, and call a path $P$ an \emph{$x,y$-path} if the vertices $x$, $y$ have degree 1 in $P$. We call $x$, $y$ the \emph{ends} of $P$. When we remove a path from a graph we will remove all the edges of the path and delete any resulting isolated vertices.

Given two disjoint vertex sets $A$ and $B$, a \emph{$d$-matching from $A$ into $B$} is a collection of disjoint sets $\{X_a\subset N(a,B):a\in A\}$ so that, for each $a\in A$, $|X_a|=d$. As is well known, such matchings can be found by showing that Hall's generalised matching condition holds. For details on this, and other standard notation, see Bollob\'as~\cite{bollo1}. We use $\log$ for the natural logarithm, and in several places omit rounding signs when they are not crucial.

\subsection{Expanders and Almost Spanning Trees}
The main properties we will use for our embeddings will be various graph expansion properties. We use the same definition of expansion as Johannsen, Krivelevich and Samotij~\cite{JKS12}.

\de Let $n\in \mathbb{N}$ and $d\in \mathbb{R}^+$. A graph $G$ is an \emph{$(n,d)$-expander} if $|V(G)|=n$ and $G$ satisfies the following two conditions.
\begin{enumerate}
\item $|N_G(X)|\geq d|X|$ for all $X\subset V(G)$ with $1\leq |X|<\lceil \frac{n}{2d}\rceil$.
\item $d_G(X,Y)>0$ for all disjoint $X,Y\subset V(G)$ with $|X|=|Y|=\lceil\frac{n}{2d}\rceil$.
\end{enumerate}
\fn

Almost spanning trees can be found in expander graphs using a theorem of Haxell~\cite{PH01}, as shown by Balogh, Csaba, Pei and Samotij~\cite{BCPS10}. We will use the following formulation of this method, which differs from that used by Johannsen, Krivelevich and Samotij~\cite{JKS12} only in that a specific vertex of the tree is embedded to a specific vertex of the expander graph. Fortunately, this version follows identically by using the full statement of the theorem of Haxell.
\begin{theorem}\label{almost}
Let $n,\Delta\in \mathbb{N}$, let $d\in \mathbb{R}^+$ with $d\geq 2\Delta$, and let $G$ be an $(n,d)$-expander. Given any tree $T\in \TT(n-4\Delta\lceil\frac{n}{2d}\rceil,\Delta)$ and vertices $v\in V(G)$ and $t\in V(T)$, we can find an embedding of $T$ in the graph $G$ with $t$ embedded on $v$.
\end{theorem}

For Theorem \ref{almostprob}, we will also require the following formulation, used by Balogh, Csaba, Pei and Samotij~\cite{BCPS10}.
\begin{theorem}\label{almost2}
Let $\Delta,m,M\in \N$. Let $H$ be a non-empty graph such that
\begin{enumerate} 
\item if $X\subset V(H)$ and $0<|X|\leq m$, then $N_H(X)|\geq \Delta|X|+1$, and 
\item if $m\leq |X|\leq 2m$, then $|N_H(X)|\geq 2\Delta|X|+M$. 
\end{enumerate}
Then $H$ contains every tree in $\TT(M,\Delta)$.
\end{theorem}

Expansion properties can also be used to construct Hamilton cycles, as shown by Hefetz, Krivelevich and Szab\'o~\cite{HKS09}.
\begin{theorem}\label{hamcycle}
Let $n,d\in \N$ satisfy that $n$ is sufficiently large and $12\leq d\leq e^{\sqrt[3]{\log n}}$. Then a graph $G$ on $n$ vertices is Hamiltonian if it satisfies the following conditions:
\begin{enumerate}
\item If $X\subset V(G)$ and $|X|\leq \frac{n\log\log n\log d}{d\log n\log\log\log n}$, then $|N(X)|\geq d|X|$,
\item If $X,Y\subset V(H)$ satisfy $|X|=|Y|\geq \frac{n\log \log n\log d}{4130\log n\log\log\log n}$ and are disjoint, then $d_G(X,Y)>0$.
\end{enumerate}
\end{theorem}
Specifically, we will use this theorem in the following form, which follows from Theorem \ref{hamcycle} by taking $d=\e\log n/\log\log n$.
\begin{corollary}\label{hamcyclesslightly}
Let $\e >0$. If $n$ is sufficiently large, and the graph $H$ is a $(n,\e\log n/\log\log n)$-expander, then $H$ is Hamiltonian.
\end{corollary}

\subsection{Finding paths}
Given many pairs of vertices in a graph and a simple large-set expansion property, we can find some path between one pair of the vertices using the following lemma.
\lem[\cite{selflove1}] \label{connect} Let $m,n\in \N$ satisfy $m\leq n/800$, let $d=n/200m$ and let $n$ be sufficiently large. Let a graph $G$ with $n$ vertices have the property that any set $A\subset V(G)$ with $|A|=m$ satifies $|N(A)|\geq (1-1/64)n$. Suppose $G$ contains disjoint vertex sets $X$, $Y$ and $U$, with $X=\{x_1,\ldots,x_{2m}\}$, $Y=\{y_1,\ldots, y_{2m}\}$ and $|U|=\lceil n/8\rceil$. Suppose, in addition, we have integers $k_i$, $i\in[2m]$, satisfying $4\log n/\log d\leq k_i\leq n/40$. Then, for some $i$, there is an $x_i,y_i$-path of length $k_i$ whose internal vertices lie in $U$.
\ma

We say a set of subgraphs \emph{covers} a graph $G$ if every vertex in $G$ is contained in one of the subgraphs. To prove Theorem \ref{comblong}, we will need the following theorem.

\begin{theorem}[\cite{selflove1}]\label{pathcoverexpander} Let $n$ be sufficiently large and let $k\in \N$ satisfy $k\geq 10^3\log^3n$, $k|n$. Let a directed graph $G$ contain $n/k$ disjoint vertex pairs $(x_i,y_i)$ and let $W=V(G)\setminus (\cup_i\{x_i,y_i\})$. Suppose $G$ has the following two properties.
\begin{enumerate}
\item  For any subset $A\subset V(G)$ with $|A|\leq n/2\log^5 n$, $|N^+(A,W)|\geq |A|\log^5 n$, and $|N^-(A,W)|\geq |A|\log^5 n$.

\item Any two disjoint subsets $A,B\subset V(G)$ with $|A|,|B|\geq n/2\log^5 n$ must have a directed edge from $A$ into $B$.
\end{enumerate}
Then we can cover $G$ with $n/k$ paths $P_i$, length $k-1$, so that, for each $i$, $P_i$ is a directed path from $x_i$ to $y_i$.
\end{theorem}

\subsection{Probabilistic results}
We will use the following results to get expansion properties in a random graph.

\begin{prop}\label{generalprops} Almost surely, if $np>20$, any two disjoint subsets $A,B\subset V(G)$ of $G=\GG(n,p)$ with $|A|=|B|=\lceil 5\log(np)/ p\rceil$ have some edge between them.
\end{prop}
\pr
Let $m=\lceil 5\log(np)/ p\rceil$. If $q$ is the probability that there exist two disjoint subsets of size $m$ which have no edge between them, then
\[
q\leq \binom{n}{m}^2(1-p)^{m^2}\leq \left(\frac{en}{m}\right)^{2m}e^{-pm^2}
\leq \left(\frac{2enp}{5\log(np)}\right)^{2m}e^{-5m\log(np) }.
\]
Therefore,
\[
q\leq (np)^{2m}e^{-5m\log(np)}= e^{-3m\log(np)}.
\]
Now, $m\log(np)\geq 5\log^2(np)/p\to \infty$ as $n\to\infty$. Therefore, $q\to 0$ as $n\to\infty$.
\oof

\begin{prop}[Alon, Krivelevich and Sudakov~\cite{AKS07}, Proposition 3.2]\label{AKS1}
Let $G=\GG(n,p)$ be a random graph with $np>20$. Then almost surely
the number of edges between any two disjoint subsets of vertices $A$, $|A|=a$, and $B$, $|B|=b$, with $abp\geq 32n$ is at least $abp/2$ and at most $3abp/2$.
\end{prop}

Our proofs of Theorems \ref{combshort} and \ref{comblong} use ideas from the sharp embedding result of Hefetz, Krivelevich, and Szab\'o~\cite{HKS12}, and we will use the following lemmas from their work.

\lem[\cite{HKS12}, Lemma 2.1] \label{HKS}Let $0<\e<1$ and $0\leq \beta\leq \e/7$ be real numbers and let $p=p(n)=(1+\e)\log n/ n$. Let $U\subset[n]$ have size $|U|\leq \b n$. Then, almost surely, the random graph $G=\GG(n,p)$ with $V(G)=[n]$ satisfies the following properties:

\begin{enumerate}
\item $\Delta(G)\leq 10\log n$.

\item $d_G(u,[n]\setminus U)\geq \eta \log n$ for every $u\in [n]$, where $0<\eta =\eta(\e)<1/2$ is a real number.
\end{enumerate} 
\ma

\lem[\cite{HKS12}, Lemma 2.4]\label{CLL} Let $G$ be a graph on $n$ vertices with maximum degree $\Delta$. Let $Y\subset V(G)$ be a set of $m=a+b$ vertices where $a$ and $b$ are positive integers. Assume that $d_G(v,Y)\geq \d$ holds for for every $v\in V$. If
\[
\Delta^2\cdot\left\lceil\frac{m}{\min\{a,b\}}\right\rceil\cdot2\cdot e^{1-\frac{\min\{a,b\}^2}{5m^2}\cdot \d}<1,
\]
then there exists a partition $Y=A\cup B$ of $Y$ such that
\begin{enumerate}
\item $|A|=a$ and $|B|=b$.
\item $d_G(v,A)\geq \frac{a}{3m}d_G(v,Y)$ for every $v\in V$.
\item$d_G(v,A)\geq \frac{b}{3m}d_G(v,Y)$ for every $v\in V$.
\end{enumerate}
\ma

We will also require the following lemma, which can be proved using a standard expectation argument, similar, for example, to calculations in the proof of Lemma \ref{HKS}.
\lem \label{mindegcond} Suppose $A\subset [n]$ and $p=p(n)$ satisfy $p|A|\geq 10\log n$. Let $d=p|A|/2$. Then almost surely the random graph $G=\GG(n,p)$ with $V(G)=[n]$ satisfies the following. For every subset $U\subset V(G)$ with $|U|\leq |A|/2d$, $|N(U,A)|\geq d|U|$.
\ma

\subsection{An important property}
The following graph property allow us to translate minimum degree conditions into expansion conditions.

\de A graph $G$ has the $(d,D,r)$-property if it contains no sets $A,B\subset V(G)$ with $|A|\leq r$, $|B|\leq d|A|$ and $d_G(A,B)\geq D|A|$.
\fn

We will typically use the $(d,D,r)$-property in the following manner. Suppose $G$ has this property and $B\subset V(G)$. Suppose further we have a set $A\subset V(G)$, with $|A|\leq r$, each of which has at least $D$ neighbours in $B$. For any subset $U\subset A$, $d(U,U\cup N(U,B))\geq D|U|$, and so, by the $(d,D,r)$-property, we must have $|N(U,B)|\geq (d-1)|U|$. That is, the subsets of $A$ expand into $B$. If, in addition, $A$ and $B$ are disjoint, then the above argument shows that $|N(U,B)|\geq |U|$ for all $U\subset A$. As Hall's generalised matching condition is satisfied, a $d$-matching from $A$ into $B$ must exist.

The following lemma, concerning when this property holds, is proved using straightforward probability, but its careful application is crucial in reaching the sharp threshold. It follows a section of the proof by Alon, Krivelevich and Sudakov of Lemma 3.1 in \cite{AKS07}.

\lem \label{mindegexp} Suppose $p=p(n)$ and $d=d(n)$ satisfies $\log^{10}n/ n\geq p\geq 4/n$ and $d\geq 4$. Let $\a,\b>0$ satisfy
\[
\a\log\left(\frac{\a}{2e\beta}\right)\geq 100.
\]
Then $G=\GG(n,p)$ almost surely has the $(d,\a d\log\log n,\max\{n/d,\beta\log\log n/ p\})$-property. 
\ma
\pr
If $G$ does not have the $(d,\a d\log\log n,\max\{n/d,\beta\log\log n/ p\})$-property, then there must exist two sets $A,B\subset V(G)$, where $|A|\leq \beta\log\log n/ p$, $|B|=d|A|$ and $d_{G}(A,B)\geq D|A|$, for $D=\a d\log\log n$ (adding vertices to $B$ if necessary to get equality). Let $p_r$ be the probability no two such sets occur with $|A|=r\leq\max\{n/d,\beta\log\log n/ p\}$. Bearing in mind that some of the edges might be counted twice if $A$ and $B$ overlap, we have
\begin{align*}\allowdisplaybreaks
 p_r &\leq \binom{n}{r}\binom{n}{dr}\binom{dr^2}{Dr/2}p^{Dr/2}
\\ &\leq\left(\frac{en}{r}\left(\frac{en}{dr}\right)^{d}\left(\frac{2edrp}{D}\right)^{D/2}\right)^r
\\ &\leq\left(\left(\frac{n}{r}\right)^{2d}\left(\frac{2edrp}{D}\right)^{D/2}\right)^r\displaybreak[4]
\\ &\leq\left(\left(\frac{2ednp}{D}\right)^{2d}\left(\frac{2edrp}{D}\right)^{D/2-2d}\right)^r
\\ &\leq\left(\log^{20d}n\left(\frac{2edrp}{D}\right)^{D/4}\right)^r.
\end{align*}
If $r<\log n$, then $2edrp/D\leq \log^{11}n/ n$ for sufficiently large $n$, and hence, as $D=\a d\log\log n$, $p_r<n^{-2}$. If $r\geq \log n$, then, as $r\leq \beta\log\log n/ p$, we have
\[
p_r\leq \log^{20dr}n\left(\frac{2e\b}{\a}\right)^{(\a dr\log\log n)/ 4}\leq \log^{(20d-100d/4)r}n\leq n^{-2}.
\]
Therefore, by looking at the sum of the probabilities $p_r$, we see the probability such a pair $A$, $B$ exists is at most $n^{-1}$.
\oof

\subsection{Dividing trees}
For Theorem \ref{combshort}, we wish to find, in a tree with lots of teeth, a much smaller subtree which still has plenty of teeth. The following is a slight generalisation of a lemma by Hefetz, Krivelevich and Szab\'o~\cite{HKS12}, and is proved below.

\lem \label{splittree}
For any $\e>0$ there exists $\b=\b(\e)>0$ and $n_0=n_0(\e)\in \N$ such that the following holds. For every tree $T$, with $|T|=n\geq n_0$, and subset $L\subset V(T)$, we can find subtrees $S, T_1, T_2\subset T$ covering $T$ so that $|S|\leq \e n$, $S$ contains at least $\b|L|$ vertices in $L$, and $T_1$, $T_2$ are disjoint and each intersects $S$ in exactly one vertex.
\ma
\begin{corollary}\label{splitspines}
For any $\e>0$ there exists $\b=\b(\e)>0$ and $n_0=n_0(\e)\in \N$ such that the following holds for any $l,k\in \N$. For every tree $T$ with $n\geq n_0$ vertices and $l$ teeth length $k$, we can find subtrees $S, T_1, T_2\subset T$ covering $T$ so that $|V(S)|\leq \e n$, $S$ has at least $(\b l-1)$ teeth length $k$ which are also teeth in $T$, and $T_1$, $T_2$ are disjoint and each intersects  $S$ in exactly one vertex.
\end{corollary}
\pr
Taking the set $L$ to be the leaves at the end of the teeth length $k$ in $T$, we apply Lemma \ref{splittree}. If the tree produced, $S$, satisfies $|V(S)\cap L|\geq 2$ then, as it is connected, it must contain all the teeth with leaves in $V(S)\cap L$.
\oof

\de
Where $S$ is a tree, we will say two subtrees $S_1$ and $S_2$ \emph{divide $S$} if they cover $S$ and intersect on precisely one vertex.
\fn
\begin{prop}\label{divide}
Given a tree $S$ we can find two trees $S_1$ and $S_2$ which divide $S$ for which $|S_1|, |S_2|\geq |S|/3$.
\end{prop}
\pr 
Take two subtrees $S_1$ and $S_2$ which divide $S$ so that $||S_1|-|S_2||$ is minimized. Let $V(S_1)\cap V(S_2)=\{v\}$.

Suppose, without loss of generality, that $|S_1|>|S_2|$. If $|S_2|\geq |S_1|-1$ then we are done, so suppose otherwise. If $v$ has only one neighbour in $S_1$, $x$ say, then the two trees on the vertex sets $V(S_1)\setminus\{v\}$ and $V(S_2)\cup\{x\}$ intersect only on $x$ and cover $S$, contradicting the choice of $S_1$ and $S_2$.

Therefore, there must be at least two neighbours of $v$ in $S_1$, and hence we can find two subtrees $S_3,S_4\subset S_1$ which cover $S_1$, each have at least two vertices and which intersect only on $v$. Without loss of generality suppose that $|S_3|\geq |S_4|$, so that $|S_4|\leq 1+(|S_1|-1)/2$. The trees on vertex sets $(V(S_1)\setminus V(S_4))\cup\{v\}$ and $V(S_2)\cup V(S_4)$ divide $S$, so, to avoid contradicting the choice of $S_1$ and $S_2$, we must have
\[
|S_1|-|S_2|\leq |S_4|-1\leq (|S_1|-1)/2.
\]
Therefore, $2|S_2|\geq |S_1|$ and so, as $|S_2|+|S_1|=|S|+1$, $|S_2|\geq |S|/3$.
\oof

\pr[Proof of Lemma \ref{splittree}] We will prove the lemma for $\e=(3/4)^{-k}$ with parameters $\b=6^{-k}$ and $n_0=12(4/3)^k$ by induction on $k\in \N\cup\{0\}$. This will prove the lemma for all $\e>0$ by taking some integer $k=k(\e)$ such that $(3/4)^{-k}<\e$. The statement holds easily for $k=0$.

Suppose then the statement holds for $k$. Given a tree $T$ with at least $12(4/3)^{k+1}$ vertices and a set $L\subset V(T)$ with $l=|L|$, find the trees $S$, $T_1$ and $T_2$ as described by the lemma for $k$ and say that $V(T_1)\cap V(S)=\{t_1\}$ and $V(T_2)\cap V(S)=\{t_2\}$. Note that if $|S|\leq (3/4)^{k+1}|T|$, then we are done, so suppose $|S|\geq (3/4)^{k+1}|T|\geq 12$. Divide $S$ into the subtrees $S_1$ and $S_2$ using Lemma \ref{divide}, so that $|S_1|,|S_2|\geq |S|/3$. Then $|S_1|,|S_2|\leq 2|S|/3+1\leq (3/4)^{k+1}|T|$. As $|S\cap L|\geq 6^{-k}|L|$, without loss of generality, we have $|S_1\cap L|\geq 6^{-k}|L|/2$. Let $V(S_1)\cap V(S_2)=\{s\}$. Note that possibly $s\in\{t_1,t_2\}$.

Suppose $t_1$, $t_2$, $s$ are distinct vertices which all lie in $S_1$. Take the common intersection vertex of the unique $t_1,t_2$-path, the $t_2,s$-path and the $t_1,s$-path in $S_1$, and call it $u$. Deleting $u$ disconnects the vertices $t_1$, $t_2$ and $s$ from each other, so we may find trees $S'_1$, $S'_2$ and $S'_3$ which intersect only on $u$, contain $t_1$, $t_2$ and $s$ respectively and cover $S$. One of these trees, $S'_j$ say, must satisfy $|S'_j\cap L|\geq 6^{-(k+1)}|L|$. Then, $S_j'$ shares at most two other vertices with the other trees $S'_i$ and the trees $S_2$, $T_1$, $T_2$, so we may proceed as below.

If the tree $S_1$ shares at most two other vertices with the trees $S_2$, $T_1$ or $T_2$ then merging any of the trees $S_2$, $T_1$ and $T_2$ which share a vertex gives at most two trees, each of which intersects with $S_1$ on precisely one vertex. Taking an additional tree consisting of a single vertex of $S_1$ if necessary, to make up a second tree, we have the induction hypothesis for $k+1$.
\oof

\section{Almost-spanning trees}\label{almostsec}

In the proof of Theorems \ref{combshort} and \ref{comblong}, we will at several points use a simple technique for finding a subgraph with expansion properties for small sets given a graph with expansion properties for large sets. This can be found in the following lemma, which we then use to prove Theorem \ref{almostprob}. 

\lem\label{neat}
Let $d\geq 1$. Suppose $G$ is a graph in which any set $A\subset V(G)$ with $|A|=m$ satisfies $|N(A)|\geq 2dm+m$. Then $G$ has a subgraph $H$, with $|H|\geq |G|-m$, in which every subset $A\subset V(G)\setminus B$ with $|A|\leq m$ satisfies $|N_H(A)|\geq d|A|$.
\ma
\pr
Let $B$ be a largest set subject to $|B|\leq m$ and $|N(B)|<d|B|$. Let $H=G[V(G)\setminus B]$. Take $A\subset V(H)$ with $0<|A|\leq m$, and suppose $|N_H(A)|<d|A|$. Then, $|N_G(A\cup B)|<d(|A|+|B|)$. By the definition of $B$, we must have that $|A\cup B|\geq m$. Thus,
\[
|N_H(A)|\geq |N_G(A\cup B)|-|N_G(B)|-|B|\geq 2dm+m-dm-m\geq d|A|.
\]
This contradicts $|N_H(A)|<d|A|$, so no such set $A$ can exist.
\oof

We will require also the following lemma from the work of Balogh, Csaba, Pei and Samotij~\cite{BCPS10}, which is proved using simple probability.

\lem[\cite{BCPS10}, Lemma 11]\label{BCPS}
Let $0<\beta\leq \gamma\leq 1/2$ and $c\geq(3/\beta)\log(4e/\gamma)$. Then almost surely the random graph $\GG(n,c/ n)$ does not contain two disjoint sets $B$ and $C$ of size at least $\beta n$ and $\gamma n$ respectively, such that $e(B,C)=0$.
\ma

\pr[Proof of Theorem \ref{almostprob}] Let $c=30(\Delta/ e)\log(1/\e)$ and $m=\lceil\e n/10\Delta\rceil$. Taking $\beta=\e/10\Delta$ and $\gamma=\e/4$ in Lemma \ref{BCPS}, the random graph $G=\GG(n,c/ n)$ almost surely has the property given in Lemma \ref{BCPS}. As there are no edges between $A\subset V(G)$ and $V(G)\setminus A$, if $|A|=m$ then $|N(A)|\geq(1-\e/4)n\geq 4\Delta m+m$.

Using Lemma \ref{neat}, find a subgraph $H\subset G$ with $|H|\geq |G|-m$ so that for any set $X\subset V(H)$, we have $|N_H(X)|\geq 2\Delta|X|\geq \Delta|X|+1$. Suppose $X\subset V(G)$ satisfies $m\leq |X|\leq 2m$. Then, by taking a subset $X'\subset X$ of size $|X'|=m$, we can see
\[
|N_{H}(X)|\geq |N_G(X')|-|X\setminus X'|-m\geq (1-\e/4)n-2m\geq (1-\e)n+4\Delta m.
\]
By Theorem \ref{almost2}, the graph $H$ contains every tree in $\TT((1-\e)n,\Delta)$, as required.
\oof

\section{$(l,\gamma)$-connectors}\label{rotate}

To construct trees with many teeth, we introduce special subgraphs called $(l,\gamma)$-connectors.

\de
A graph $H$ is an \emph{$(l,\gamma)$-connector} if $H$ has $l$ vertices and there are two disjoint subsets $H^+,H^-\subset V(H)$ with size $|H^+|=|H^-|= \lceil\gamma l\rceil$ so that given any pair of vertices $x\in H^+$ and $y\in H^-$ there is an $x,y$-Hamilton path in $H$. 
\fn

In practice, when we have found an $(l,\gamma)$-connector, $H$ say, we will implicitly fix two such sets $H^+$ and $H^-$. We will often treat a $(l,\gamma)$-connector $H$ as if it were a normal path with length $l-1$ and ends $H^+$ and $H^-$. We can then connect two vertices $x$ and $y$ through $H$ if we can find an edge from $x$ to $H^+$ and an edge from $y$ to $H^-$ by taking the associated Hamilton path through $H$. The advantage of using $(l,\gamma)$-connectors, as we will see, is that we are more likely to be able to find these edges than if we were using just an ordinary path with two end vertices. We will find the following definition useful in describing such connections.

\de
Given a collection of subsets $\mathcal{A}$ in the graph $G$, the \emph{grouped graph on the set $\mathcal{A}$ with respect to $G$} is the graph $H$ with vertex set $\mathcal{A}$ and an edge between $B,C\in \mathcal{A}$ if there is some edge between $B$ and $C$ in $G$.
\fn
We will often deal with bipartite grouped graphs with vertex classes $\mathcal{A}$ and $\mathcal{B}$, where we only consider edges between the classes. In a slight abuse of notation, when we say a bipartite grouped graph has vertex classes $\mathcal{A}$ and $V$ where $\mathcal{A}$ is a collection of subsets and $V\subset V(G)$, we shall mean such a graph on vertex classes $\mathcal{A}$ and $\{\{v\}:v\in V\}$. The following notation will also be useful.
\de
Given a collection $\mathcal{A}$ of $(l,\gamma)$-connectors, let $\mathcal{A}^+=\{P^+:P\in \mathcal{A}\}$ and $\mathcal{A}^-=\{P^-:P\in \mathcal{A}\}$.
\fn

An $(l,\gamma)$-connector can be found in any graph with a simple expansion property, as follows.
\lem \label{rotatingpaths}
Let $\lambda>0$ and let $n$ be sufficiently large, based on $\lambda$. Let $G$ be a graph on $n$ vertices in which any set $A\subset V(G)$ with $|A|=\lceil\lambda n/\log n\rceil$ satisfies $N(A)\geq (1-1/64)n$. If $l_0\leq n/6$, then $G$ contains an $(l_0,\gamma)$-connector with $\gamma=\log\log n/ 16\log n$.
\ma

Before proving Lemma \ref{rotatingpaths}, we will sketch an almost sure embedding of $\mathrm{Comb}_{n,k}$ in the random graph $\GG(n,\log^2 n/n)$, when $k|n$. We then give an indication of how we will reduce the probability required for the embedding. We will take $k=\sqrt{n}$, but the method works for general $k$.

\subsection{Embedding $\mathrm{Comb}_{n,\sqrt{n}}$ in $\GG(n,\log^2n/ n)$}
\lem\label{log2sketch} If $\sqrt{n}\in \N$, then almost surely there is a copy of $\mathrm{Comb}_{n,\sqrt{n}}$ in $G=\GG(n,\log^2n/ n)$
\ma
\pr[Sketch proof of Lemma \ref{log2sketch}]
We reveal the edges of $G$ in three rounds, where at each stage any edge is present independently with probability $p=\log^2 n/3n$. The edges of the final graph appear then with probability at most $1-(1-p)^2\leq 3p$. Therefore, if a copy of $\mathrm{Comb}_{n,\sqrt{n}}$ exists almost surely in such a random graph, then almost surely one must exist in the random graph $\GG(n,\log^2 n /n)$.

Reveal the first set of edges to get $G_1$. Almost surely, any subgraph $H\subset V(G)$ with $|H|\geq n/3$ will satisfy the conditions of Lemma \ref{rotatingpaths} with $\lambda=1$ (using, for example, Proposition \ref{AKS1}). Therefore, any subset of $n/3$ vertices must contain an $(l,\gamma)$-connector for $\gamma=\log\log n/16\log n$ and $l=\lceil(\sqrt{n}-1)/2\rceil$. We may then find greedily in the graph a set $\mathcal{A}$ of $\sqrt{n}$ disjoint $(l,\gamma)$-connectors. Find a path length $\sqrt{n}-1$ and label it $Q=q_1\ldots q_{\sqrt{n}}$ (for example by taking a path through a $(\sqrt{n},\gamma)$-connector). This path will form the `spine' of our comb.

Let $W$ be the set of vertices in the graph not in any of the connectors or the path $Q$, so that $|W|=\sqrt{n}\lfloor(\sqrt{n}-1)/2\rfloor\geq n/3$. By revealing more edges with probability $p$ to get $G_2$ we can almost surely find a Hamilton cycle in $G_2[W]$. Take this cycle and break it into $\sqrt{n}$ paths length $\lfloor(\sqrt{n}-1)/2\rfloor-1$. Label these paths by $R_i$, $i\in[\sqrt{n}]$, and in each path pick an end vertex and label it $r_i$.

We have now covered our graph by a path of length $\sqrt{n}-1$, $Q$, the $\sqrt{n}$ $(\lceil(\sqrt{n}-1)/2\rceil,\gamma)$-connectors in $\mathcal{A}$, and $\sqrt{n}$ paths of length $\lfloor(\sqrt{n}-1)/2\rfloor-1$, $R_i$, as illustrated in Figure \ref{absorbpic}. Suppose there is a matching in the bipartite grouped graph on vertex sets $\{q_i:i\in[\sqrt{n}]\}$ and $\mathcal{A}^+$, and a matching in the bipartite grouped graph on vertex sets $\mathcal{A}^-$ and $\{r_i:i\in[\sqrt{n}]\}$. Then, for each $i$, we could take the connector $P\in \mathcal{A}$ for which $P^+$ is matched to $q_i$, and take the vertex $r_j$ matched to $P^-$, and find a neighbour $q$ of $q_i$ in $P^+$ and a neighbour $r$ of $r_j$ in $P^-$. Taking a $q,r$-Hamilton path in $G_1[P]$, we could then attach a tooth length $\sqrt{n}-1$ to $q_i$ using this path and $R_j$ along with the edges $q_iq$ and $rr_i$. Doing this for each $i\in[\sqrt{n}]$ gives a copy of $\mathrm{Comb}_{n,\sqrt{n}}$. 

These matchings almost surely exist if we reveal more edges with probability $p$. Indeed, the probability an edge in the two grouped graphs is present is at least
\[
1-(1-p)^{\gamma l}\geq p\gamma l/2\geq \frac{\log n\log\log n}{16\sqrt{n}}.
\]
Thus the two grouped graphs are random bipartite graphs with equal class sizes $m=\sqrt{n}$ and edges present independently with probability at least $2\log m /m$. This is above the threshold for a matching to almost surely exist in such a graph (see, for example, Bollob\'{a}s~\cite{bollo1}).
\oof

\setlength{\unitlength}{0.05cm}
\begin{figure}[t]
\centering
\scalebox{1.2}{
\begin{picture}(100,120)

\linethickness{0.04cm}
\multiput(10,50)(20,0){5}{\circle*{5}}
\multiput(10,105)(20,0){5}{\circle*{5}}
\multiput(10,15)(20,0){5}{\line(0,1){35}}

\multiput(5,60)(20,0){5}{\line(1,0){10}}
\multiput(5,70)(20,0){5}{\line(1,0){10}}
\multiput(5,60)(20,0){5}{\line(0,1){10}}
\multiput(15,60)(20,0){5}{\line(0,1){10}}

\multiput(5,85)(20,0){5}{\line(1,0){10}}
\multiput(5,95)(20,0){5}{\line(1,0){10}}
\multiput(5,85)(20,0){5}{\line(0,1){10}}
\multiput(15,85)(20,0){5}{\line(0,1){10}}

\multiput(10,70)(20,0){5}{\line(0,1){15}}

\put(10,105){\line(1,0){80}}

%\put(187,21){$y_0$}

%\put(140,30){\line(1,0){40}}

\put(-6,63){$P_1^-$}
\put(-6,88){$P_1^+$}
\put(48,110){$Q$}
\put(8,110){$q_1$}
\put(-2,75.5){$P_1$}
\put(0,48){$r_1$}
\put(-1,29.5){$R_1$}

\linethickness{0.01cm}
\qbezier(10, 105)(18,100)(26, 95)
\qbezier(30, 105)(38,100)(46, 95)
\qbezier(50, 105)(58,100)(66, 95)
\qbezier(70, 105)(78,100)(86, 95)

\qbezier(15, 93)(32.5,99)(50, 105)
\qbezier(35, 93)(52.5,99)(70, 105)
\qbezier(55, 93)(72.5,99)(90, 105)

\qbezier(65, 93)(47.5,99)(30, 105)
\qbezier(10, 95)(10,100)(10, 105)
\qbezier(10, 95)(10,100)(10, 105)

\qbezier(30, 95)(30,100)(30, 105)

\qbezier(10, 50)(18,55)(26, 60)
\qbezier(30, 50)(38,55)(46, 60)
\qbezier(50, 50)(58,55)(66, 60)
\qbezier(50, 50)(42,55)(34, 60)
\qbezier(30, 50)(22,55)(14, 60)
\qbezier(70, 50)(78,55)(86, 60)

\qbezier(15, 62)(32.5,56)(50, 50)
\qbezier(35, 62)(52.5,56)(70, 50)
\qbezier(55, 62)(72.5,56)(90, 50)

\qbezier(10, 50)(10,55)(10, 60)

%\qbezier(65, 93)(47.5,99)(30, 105)
%\qbezier(10, 95)(10,100)(10, 105)
%\qbezier(10, 95)(10,100)(10, 105)

%\qbezier(140, 30)(185, 58)(190, 30)%165%60

%\put(20,30){\line(1,0){170}}
\end{picture}
}
\vspace{-0.75cm}
\caption{An almost sure embedding of $\mathrm{Comb}_{n,\sqrt{n}}$ in $\GG(n,\log^2 n/n)$.}
\label{absorbpic}

\end{figure}

%\begin{figure}[h!]\label{comb}
 % \centering
  %    \includegraphics[width=6cm]{sketchlog2}

%\end{figure}

The limiting requirement for the probability used in this sketch is in finding the matchings, which is limited by the need to ensure that every vertex $q_i$ or $r_i$ has some neighbour in a set $P^+$ or $P^-$ respectively for some connector $P\in \mathcal{A}$. If the probability $p=C\log n/ n$, for any constant $C$, is used then we must expect vertices which have no neighbours in the sets $P^+$ or $P^-$ for any $P\in \mathcal{A}$. However, on average vertices will still have many such neighbours, at least $c\log\log n$, for some small $c$. If we can construct the paths $Q$ and $R_i$ so that all the vertices $r_i$ and $q_i$ have above the average of such neighbours then, using Lemma \ref{mindegexp}, finding the matchings will be easier. The proof of Theorem \ref{combshort} is more involved, but this is the basic idea behind reducing the probability required.

\subsection{Constructing $(l,\gamma)$-connectors}
Our construction of $(l,\gamma)$-connectors is inspired by the celebrated technique of P\'osa rotation, introduced by P\'osa to find Hamilton cycles in random graphs~\cite{posa76}. Though we do not use P\'{o}sa rotations explicitly, essentially we construct our $(l,\gamma)$-connector as a path with additional edges that guarantee the path can be rotated at both ends, independently, to give at least $\gamma l$ new end vertices.

\pr[Proof of Lemma \ref{rotatingpaths}]
Let $k=4\log n / \log\log n$, $m=\lfloor\lambda n/\log n\rfloor$, and $d=n/200m\geq \log n/200\lambda$. Assume that $l_0\geq 2k$, otherwise we may simply take a path length $l_0$ (found using, for example, Lemma \ref{connect}) as an $(l_0,\gamma)$-connector.

We will call a path $P$ starting at $p_0$ a \emph{good path} if there is a set $R\subset E(P)$ such that
\begin{itemize}
\item $|R|\geq |P|/ 2k$, and,
\item for each $e\in R$, there is a Hamilton path in $G[V(P)]$ which starts at $p_0$, passes through each edge in $R$ and whose last edge is $e$.
\end{itemize}

We say that such a set $R$ \emph{demonstrates that $P$ is a good path}. Let $S$ be any set in $G$ containing at least $n/3$ vertices. 

\begin{claim}
Suppose we have a collection of disjoint good paths $P_i$, $i\in I$, in $G[S]$ with initial vertices $p_i$, such that
\begin{equation}\label{pathprop}
2m\leq \sum_i\left\lceil \frac{|P_i|}{2k}\right\rceil\quad\quad\text{and}\quad\quad\sum_i|P_i|\leq \frac{n}{6}.
\end{equation}
Suppose also we have a set of integers $k_i$, $i\in I$, with $k-1\leq k_i\leq 2k-1$. Then
we can find a new path $P'_j$, for some $j\in I$, with initial vertex $p_j$, which is a good path in $S$, is disjoint from all the paths $P_i$, $i\neq j$, and contains $k_j$ more vertices than $P_j$ does.
\end{claim}
In other words, given such a collection of paths satisfying (\ref{pathprop}) we can lengthen one of the paths using vertices from $S$ so that it is still a good path with the same starting vertex, keeping the paths disjoint.

\pr[Proof of the claim]
Suppose we have such a collection of disjoint good paths $P_i$, $i\in I$. For each path $P_i$, find a set $R_i$ which demonstrates that it is a good path. By (\ref{pathprop}), $\sum_i|R_i|\geq \sum_i\lceil |P_i|/2k\rceil \geq 2m$. Let $S'=S\setminus(\cup_iV(P_i))$, $R=\cup_iR_i$ and $r=|R|$. Label the vertices which appear in the edges in $R$, so that $R=\{x_iy_i:1\leq i\leq r\}$ and so that if $x_iy_i\in R_j$ then $x_i$ appears earlier in the path $P_j$ starting from $p_j$ than $y_i$ does. This labelling ensures that the vertices $x_i$ are distinct and that the vertices $y_i$ are also distinct. For each $i$, find $j$ for which $x_iy_i\in R_j$ and let $k'_i=k_j$.

By Lemma \ref{connect}, as $|R|\geq 2m$ and $|S'|\geq n/6$, for some $i$ there is an $x_i$,$y_i$-path $Q$ of length $k'_i+1$ whose interior vertices lie in $S'$. Say that the edge $x_iy_i$ is in the path $P_j$, so that $Q$ has length $k'_i+1=k_j+1$. 

Let $P_j'$ be the path formed by replacing the edge $x_iy_i$ in $P_j$ by the path $Q$. This lengthens $P_j$ by the correct number of vertices, and,  to finish the proof of the claim, we need only show that $P_j'$ is a good path. We will do this by finding a demonstrating set.

As $P_j$ is a good path, there is a path $P'$ through $V(P_j)$ starting at $p_j$, which passes through each edge in $R_j$ and whose final edge is $x_iy_i$. Switch the labelling of $x_i$ and $y_i$, if necessary, so that this path ends in $y_i$. Label the vertices of $Q$ so that $Q$ is the path $x_iq_1q_2\ldots q_{k_j}y_i$.

Let $R'=(R_j\setminus \{x_iy_i\})\cup\{q_1q_2,q_{k_j}y_i\}$. Then $|R'|=|R_j|+1\geq (|V(P_j)|+k_j)/2k$. Let $e\in R_j\setminus\{x_iy_i\}$. By the definition of $R_j$, there is a Hamilton path in $G[V(P_j)]$ which passes through all the edges of $R_j$ and ends in $e$. This path must contain the edge $x_iy_i$, so, by replacing that edge with the path $Q$, we get a Hamilton path in $G[V(P'_j)]$ which passes through all the edges of $R'$ and ends with $e$. Recalling the path $P'$, replace the edge $x_iy_i$ by $Q$ to get a Hamilton path in $G[V(P'_j)]$ which passes through all the edges in $R'$ and ends in $q_{k_j}y_i$. If we add instead the path $y_iq_{k_j}\ldots q_2q_1$ to the end of $P'$ then we get a Hamilton path in $G[V(P'_j)]$ which passes through all the edges in $R'$ and ends with $q_1q_2$. Thus, $R'$ demonstrates that $P'_j$ is a good path, and gives the claim.
\oof

Let $l_1=\lfloor (l_0-1)/2\rfloor$. Given any collection of $2m$ edges $\{x_iy_i:i\in I\}$ and a set $S$ of at least $n/3$ vertices disjoint from these edges, we claim we can find a good path of length $l_1$, which has $x_j$ as its initial vertex, for some $j$, and which lies in $S\cup\{x_j,y_j\}$.

To show this, consider each edge $x_iy_i$ as a path $P_i$ of length 1, where the set $\{x_iy_i\}$ demonstrates that this is a good path. Thus, the edges $x_iy_i$ satisfy the conditions of the claim above. 

Repeatedly apply the claim with $k_i=k-1$ if $l_1\geq |P_i|+2k-1$, and $k_i=l_1-|P_i|$ otherwise. The paths lengthen until either one of them has length $l_1$, or until the upper bound in (\ref{pathprop}) is not satisfied. In the latter case, we discard a shortest path repeatedly until the upperbound in (\ref{pathprop}) holds again. By repeatedly applying the claim, and then removing paths, eventually we must find a path length $l_1$, as required.

This claim also holds identically with $l_2=\lceil l_0/2\rceil$ in place of $l_1$.

We can now build our $(l_0,\gamma)$-connector. First, we find in the graph $G$ $4m$ disjoint paths with length 2 and label their vertices $y_ix_iy'_i$, $i\in[4m]$. If $U\subset V(G)$ satisfies $|U|\geq n/3$, then pick two disjoint subsets $U_1,U_2\subset U$ with $|U_1|=|U_2|=m$. As $|N(U_1,U)\cap N(U_2,U)|\geq |U|-n/32 -2m$, using the expansion property for $G$, we can certainly pick vertices $v\in U$, $u_1\in U_1$ and $u_2\in U_2$ so that $vu_1$, $vu_2\in E(G)$. Therefore in any subset of $G$ of size at least $n/3$ we can find a path with length 2, so we may greedily select the paths described.

Divide the vertices not in these short paths into two sets $S_1$ and $S_2$ of size at least $n/3$. Given any subset $M\subset[4m]$ of size $2m$ we can find an index $i\in M$ and a good path $P$ in $S_1\cup\{x_i,y_i\}$ with length $\lfloor (l_0-1)/2\rfloor$ which starts with $x_i$. Therefore this must be possible for at least $2m+1$ values of $i\in[4m]$. Similarly, for at least $2m+1$ values of $i$ there must be a good path from $x_i$ in $S_2\cup\{x_i,y'_i\}$ with length $\lceil (l_0-1)/2 \rceil$. There must be then some index $j\in[4m]$ and paths $P_1$ and $P_2$ which are good, start on $x_j$, are disjoint except for the vertex $x_j$, and have length $\lceil (l_0-1)/2\rceil$ and $\lfloor (l_0-1)/2\rfloor$ respectively. Let $H$ be the subgraph $G[V(P_1)\cup V(P_2)]$.

Let $R_1$ and $R_2$ be sets demonstrating that $P_1$ and $P_2$ respectively are good paths. For each $e\in R_1$ there is vertex $v$ in $e$ for which there is a Hamilton path in $G[V(P_1)]$ from $x_j$ to $v$ which goes through every edge in $R_1$. Pick such a vertex and call it $v_e$. Let $H^+=\{v_e:e\in R_1\}$. Note that if $v_e=v_{e'}$ then there is a path going through all the edges in $R_1$ which ends in $v_e$, so $e=e'$. Therefore $|H^+|\geq (l_0-1)/4k\geq \gamma l_0$. Define similarly vertices $v_e$ for each edge $e\in R_2$ and let $H^-=\{v_e:e\in R_2\}$, so that $|H^-|\geq \gamma l_0$.

For any pair of vertices $x\in H^+$ and $y\in H^-$ we may combine an $x,x_j$-Hamilton path in $G[V(P_1)]$ with a $x_j,y$-Hamilton path in $G[V(P_2)]$ and get a path from $x$ to $y$ covering exactly the vertices in $H$. Thus, $H$ is an $(l_0,\gamma)$-connector.
\oof

\section{Proof of Theorem \ref{combshort}}\label{secshort}

Before proving Theorem \ref{combshort}, we will prove the following useful lemma. Given several graph properties, Lemma \ref{embedavoid} embeds two trees into a graph before covering the rest of the vertices with paths of the same length, so that a set of `bad' vertices are contained within these paths.

\lem\label{embedavoid} Let $\e,\eta>0$ be fixed, and let $n\in\N$ be sufficiently large, depending on $\e$ and $\eta$. Let $l\geq 5$ and $s\in \N$ satisfy $ls\geq \e n$. Let $T_1$ and $T_2$ be trees so that $|T_1|+|T_2|=n-sl$, and suppose they have vertices $t_1\in V(T_1)$, and $t_2\in V(T_2)$. Suppose a graph $G$ has $n$ vertices and contains the set $Z\subset V(G)$ with $|Z|\geq n-l s/8$, so that the following holds.
\begin{enumerate}
\item For every vertex $v\in V(G)$, $d(v,Z)\geq \eta\log n$,

\item $G$ has the $(d,D,r)$-property for $d=\eta\e\log n/10^5\log\log n$, $D=\eta\e\log n/ 3^6$ and $r=10n\log\log n/ \log n$,

\item Every two disjoint subsets of $V(G)$ of size $10n\log\log n/\log n$ have some edge between them, and,

\item $\Delta(G)\leq 10\log n$.

\end{enumerate}
Then, given any two distinct vertices $v_1,v_2\in V(G)\setminus Z$, we may embed the trees $T_1$ and $T_2$ disjointly in $Z\cup\{v_1,v_2\}$, so that $t_1$ and $t_2$ are embedded on $v_1$ and $v_2$ respectively, and so that the vertices not used in the embedding of $T_1$ or $T_2$ span $s$ paths of length $l-1$ with end vertices in $Z$.
\ma
\pr
Let $B=V(G)\setminus Z$. The conditions for Lemma \ref{CLL} comfortably hold for $G$ and $Z$, as $e^{\eta\log n}$ is much larger than $\Delta(G)$. By applying Lemma \ref{CLL} three times, we can therefore partition $Z$ into the sets $Z_1,\ldots, Z_5$, as follows. We have $|Z_1|=2n/3$, $|Z_3|= \lfloor sl/8\rfloor$, $|Z_4|=\lfloor s/2\rfloor\lfloor l/2 \rfloor+2\lfloor s/2\rfloor$, $|Z_5|=\lceil s/2\rceil \lceil l/2\rceil$, with $|Z_2|$ determined by the partitioning of $Z$, and for each $v\in V(G)$ and each $i$,
\[
d(v,Z_i)\geq \frac{|Z_i|}{3^4|Z_0|}\eta \log n\geq D,
\]
where we have used that, for each $i$, $|Z_i|\geq sl/9\geq \e n/9$.

Suppose $U\subset V(G)$ with $|U|\leq r$. Each vertex $u\in U$ has $d(u,Z_i)\geq D$, for each $i$, so, as $G$ has the $(d,D,r)$-property, we must have $|N(U,Z_i)|\geq (d-1)|U|$. Disjoint sets of size $r$ have some edge between them, and, for each $i$, $r\leq |Z_i|/2d$. Therefore, if a vertex set $W$ contains $Z_i$, then $G[W]$ is a $(|W|,d-1)$-expander. For each $i$, and any set $U\subset V(G)\setminus Z_i$ with $|U|\geq m$, we have $|N(U,Z_i)|\geq |Z_i|-m\geq |Z_i|/2$. Therefore, using Hall's matching condition, for each $i$, any set $U\subset V(G)\setminus Z_i$ with $|U|\leq |Z_i|/2$ must have a matching into $Z_i$.

Suppose, without loss of generality, that $T_1$ is smaller than $T_2$, so that it contains at most $n/2$ vertices. Let $W_1=Z_1\cup\{v_1\}$. As $W_1$ is a $(|W_1|,d-1)$-expander and $v_1\in W_1$, by Theorem \ref{almost} there is a copy of $T_1$, say $S_1$, in $G[W_1]$ in which $t_1$ is embedded on $v_1$. Let $W_2=(W_1\setminus V(S_1))\cup Z_2\cup\{v_2\}$. The graph $G[W_2]$ is a $(|W_2|,d-1)$-expander because it contains $Z_2$. It also contains the vertex $v_2$, and in total $|W_2|=|T_2|+sl-|Z_3\cup Z_4\cup Z_5\cup B|+1\geq |T_2|+sl/8\geq |T_2|+\e n/8$ vertices. Therefore, by Theorem \ref{almost}, $G[W_2]$ contains a copy of $T_2$ with $t_2$ embedded to $v_2$. Let such a copy of $T_2$ be $S_2$.

Let $W_3=V(G)\setminus (Z_4\cup Z_5\cup V(S_1)\cup V(S_2))$, so that $(B\setminus\{v_1,v_2\})\cup Z_3\subset W_3$. Now, $|W_3|=n-|S_1|-|S_2|-\lceil s/2\rceil l-2\lfloor s/2\rfloor=\lfloor s/2\rfloor(l-2)$. By Corollary \ref{hamcyclesslightly}, as $G[W_3]$ is a $(|W_3|,d-1)$-expander it must contain a Hamilton cycle, if $n$ is sufficiently large. Take such a Hamilton cycle and break it into $\lfloor s/2\rfloor$ paths length $l-3$. These paths cover the vertices in $Z_3$, so contain the vertices in $B\setminus\{v_1,v_2\}$. Take the end vertices of all these paths and find a matching from them into $Z_4$. The paths now have length $l-1$ and ends lying in $Z_4$, which is disjoint from $B$.

Finally, let $W_4$ be the vertices from $Z_4$ not used as good end vertices for the paths and add the vertices from $Z_5$, so that $|W_4|=\lceil s/2\rceil l$. As $Z_5\subset W_4$, $G[W_4]$ is a $(|W_4|,d)$-expander, so we may find a Hamilton cycle in $G[Z_5]$ with Corollary \ref{hamcyclesslightly}. Break this Hamilton cycle into paths of length $l-1$ to complete the structure required in the lemma.
\oof

Theorem \ref{combshort} is proved in two different cases. The case when $k\geq \log^9 n$ is covered by Theorem \ref{comblong}, so we will focus on proving Theorem \ref{combshort} when $k\leq \log^9 n$. However it is possible to prove Theorem \ref{combshort} in full without using the tools used for Theorem \ref{comblong} and we will sketch how this can be done at the end of this section. We will also presume that $k\geq 16\log n/\log\log n$, before remarking the small change necessary to cover the case when $10\leq k\leq 16\log n/\log\log n$.

\pr[Proof of Theorem \ref{combshort} when $16\log n/\log\log n\leq k\leq \log^9 n$] 
Let $T$ be a tree with $n$ vertices and at least $\a n/k$ teeth length $k$.

We will expose the edges of $G$ in five stages to get the graphs $G_1$, \ldots, $G_5$. In the first four stages we will reveal edges with probability $\e\log n/8n$, and in the final stage we will reveal edges with probability $(1+\e)\log n/2n$. Therefore, in total, the probability any edge is present in the final graph is
\[
1-\left(1-\frac{\e\log n}{8n}\right)^4\left(1-\frac{(1+\e)\log n}{2n}\right)\leq \frac{(1+\e)\log n}{n}.
\]
Thus, if we can almost surely construct a copy of $T$ in the resulting graph $G$, then almost surely we can find a copy of $T$ in the random graph $\GG(n,(1+\e)\log n/n)$.

Let $\beta=\beta(\e/14)$ come from Lemma \ref{splittree}. Let $\mu_1\leq \beta\a/4$ be sufficiently small that, by Lemma \ref{mindegexp}, a random graph $\GG(n,p)$, for any $p=p(n)\leq \log^{10}n/ n$, almost surely has the $(4,\e\log\log n /10^4,\mu_1\e\log\log n / p)$-property. Let $\mu_2\leq\mu_1/20$ be sufficiently small that, by Lemma \ref{mindegexp}, a random graph $\GG(n,p)$, for any $p\leq \log^{10}n/ n$, almost surely has the ($4$,$\e\mu_1\log\log n /800$, $\mu_2\e\log\log n /p$)-property.

Reveal edges with probability $\e \log n /8n$ to get the graph $G_1$. Let $m_1=\lceil 50n\log\log n /\e\log n \rceil$. By Proposition \ref{generalprops}, almost surely any two disjoint sets $A,B\subset V(G)$ with $|A|,|B|\geq m_1$ must have $d_{G_1}(A,B)>0$. By Proposition \ref{AKS1}, almost surely every set $A\subset V(G)$ with $|A|=10^6 n/\e\log n$ satisfies $|N_{G_1}(A)|\geq 127n/128$. Therefore, for sufficiently large $n$, any vertex set of size $n/2$ satisfies the conditions of Lemma \ref{rotatingpaths} with $\lambda=\e/10^5$ and so contains an $(l,\gamma)$-connector, where $\gamma=\log\log n/ 16\log n$ and $l=\lfloor k/2\rfloor$. Greedily then, we may find two sets $\mathcal{A}$ and $\mathcal{B}$, each of which contains $\lfloor n/4l\rfloor$ $(l,\gamma)$-connectors in $G_1$, where these connectors are all disjoint.

Reveal edges with probability $\e\log n /8n$ to get $G_2$. Let $H$ be the bipartite grouped graph on $\mathcal{A}^-$ and $\mathcal{B}^+$. Edges are present in $H$ independently with probability $p_1$, where
\[
p_1=1-\left(1-\frac{\e\log n}{8n}\right)^{\left\lceil \frac{l\log\log n}{16\log n}\right \rceil^2}\geq \frac{\e\log n}{16n}\left\lceil \frac{l\log\log n}{16\log n} \right\rceil^2\geq \frac{\e l \log\log n}{256 n}.
\]
Let $m_2=\lceil\mu_1n/10^3l\rceil$. By Proposition \ref{generalprops}, almost surely we can assume that, given any two subsets $\mathcal{Y}\subset \mathcal{A}$ and $\mathcal{Z}\subset \mathcal{B}$ with $|\mathcal{Y}|, |\mathcal{Z}|\geq m_2$, we have $d_H(\mathcal{Y}^-,\mathcal{Z}^+)>0$. Greedily select $|\mathcal{A}|/2$ independent edges in the graph $H$. This is possible, otherwise removing the vertices associated with a maximal set of independent edges would leave a graph with no edges and at least $|\mathcal{A}|/2$ remaining vertices in each half, a contradiction. Let $\mathcal{A}_1\subset \mathcal{A}$ be a set of connectors $P\in \mathcal{A}$ for which $P^-$ appears in the independent edges, with $|\mathcal{A}_1|=\lfloor n/20l\rfloor$.

Reveal edges with probability $\e\log n/ 8n$ to get the graph $G_3$. Let $K$ be the bipartite grouped graph with vertex classes $\AA_1^+$ and $V:=V(G)\setminus (\cup_{P\in \AA\cup \BB}V(P))$. This graph has edges present independently with some probability $p_2$, where
\begin{equation}\label{bounds}
\frac{\e l\log\log n}{400n}\leq p_2=\left(1-\frac{\e\log n}{8n}\right)^{\left\lceil \frac{l\log\log n}{16\log n}\right \rceil}\leq \frac{\e l\log\log n}{n}.
\end{equation}
The upper bound in (\ref{bounds}) relies on $\lceil l\log\log n/ 16\log n\rceil$ being at most a factor of $2$ greater than $l\log\log n/ 16\log n$, which holds as $l\geq 16\log n /\log\log n$. 

Let $D_1=\e\log\log n/10^4< p_2|\AA_1|/2$. As $\mu_1\e\log\log n / p_2\geq \mu_1n/ l$, by Lemma \ref{mindegexp}, the choice of $\mu_1$ and considering the graph $K$ as a subgraph of the random graph $\GG(n,p_2)$, $K$ almost surely has the $(4,D_1,\mu_1n/l)$-property.

Let $V_1$ be the set of vertices in $V$ which have at least $D_1$ neighbours in $\mathcal{A}_1$ in the graph $K$. There are at most $D_1|V\setminus V_1|< p_2|\mathcal{A}_1||V\setminus V_1|/2$ edges between $\mathcal{A}_1$ and $V\setminus V_1$ in the graph $K$. Almostly surely then, by Proposition \ref{AKS1}, $2D_1|V\setminus V_1|\leq 32 n$, and hence $|V\setminus V_1|=o(n)$.

The set $V$ contains all the vertices not in any of the connectors, so that $|V|\geq n/2$. Recall that there is an edge between any two disjoint vertex sets of size $m_1$ in $G_1$, and let $d_1=n/100m_1\geq\e\log n/10^4\log\log n$. Given any set $U\subset V$ of size $m_1$ then, $|N_{G_1}(U,V)|\geq |V|-2m_1\geq 2d_1m_1+2m_1$. Therefore, by Lemma \ref{neat}, there is some subset $V_2\subset V_1$ of size at least $|V_1|-m_1\geq n/3$, so that given any set $U\subset V_2$, with $|U|\leq m_1$, $|N(U,V_2)|\geq d_1|U|$. Thus, $G_1[V_2]$ is a $(|V_2|,d_1)$-expander.

We will now start our embedding. Using Corollary \ref{splitspines}, split $T$ into three trees $S$, $T_1$, and $T_2$ covering $T$ so that $T_1$ and $T_2$ are disjoint, $V(S)\cap V(T_1)=\{t_1\}$, $V(S)\cap V(T_2)=\{t_2\}$, and $S$ contains at least $\a\b n/2k$ teeth length $k$ which are also teeth in $T$ (i.e.\ they do not contain $t_1$ or $t_2$), but at most $\e n/14$ vertices.

Pick $s_1:=\lfloor\mu_1n/l\rfloor$ teeth in $S$ with length $2l\leq k$, possible as $\mu_1\leq \a\b/4$, and remove them to leave the tree $S'$. This tree has maximum degree at most $\Delta$ and at most $\e n/14$ vertices, so we may embed $S'$ in $G[V_2]$ using Theorem \ref{almost}. Let $X$ be the set of $s_1$ vertices in $V(G)$ to which we need to attach teeth length $2l$ to extend this embedding to $S$.

As each vertex in $X$ lies in $V_1$, it must have at least $D_1$ neighbours in $\mathcal{A}^+$ in the graph $K$. Therefore, as $|X|=s_1$ and $K$ has the $(4,D_1,s_1)$-property, we can find a matching from $X$ into $\mathcal{A}^+$ in $K$. This matching allows us to attach a different $(l,\gamma)$-connector $P\in \AA$ to each vertex in $X$ using the set $P^+$. Let $\mathcal{A}_2\subset \mathcal{A}_1$ be the set of connectors attached in this manner. We have sucessfully embedded the tree $S$ except for $s_1$ teeth length $2l$ which are replaced by $(l,\gamma)$-connectors.

Take a subset $\mathcal{B}_3\subset \mathcal{B}_2$, so that $|\mathcal{B}_3|=|\mathcal{A}_2|/2$. For every set $\mathcal{U}\subset \mathcal{A}_2\cup\mathcal{B}_3$ with $|\mathcal{U}|=2m_2$, either $|\mathcal{U}\cap\mathcal{B}_3|\geq m_2$ or $|\mathcal{U}\cap \mathcal{A}_2|\geq m_2$. Therefore we have $|N_H(\mathcal{U},\mathcal{A}_2\cup\mathcal{B}_3)|\geq \min\{|\mathcal{A}_2|-2m_2,|\mathcal{B}_3|-2m_2\}\geq |\mathcal{B}_3|-2m_2\geq 20m_2$. By Lemma \ref{neat}, there is a subset $\mathcal{Y}\subset\mathcal{A}_2\cup\mathcal{B}_3$ with $|\mathcal{Y}|\geq|\mathcal{A}_2\cup\mathcal{B}_3|-2m_2$, so that, for all $\mathcal{U}\subset \mathcal{Y}$ with $|\mathcal{U}|\leq 2m_2$ we have $N_H(\mathcal{U},\mathcal{Y})\geq 4|\mathcal{U}|$. Let $\AA_3=\mathcal{Y}\cap\AA_2$, and $\BB_4=\mathcal{Y}\cap \BB_3$.

For each connector $P\in \AA_2\setminus \AA_3$, use the independent edges between $\AA_2^-$ and $\BB_1^+$ to attach an $(l,\gamma)$-connector to $P$. If $P'$ is attached to $P$ in this manner, then by taking the vertex $x\in X$ attached to $P$, we can find a path length $2l$ starting from $x$ and covering the two connectors $P$ and $P'$. This attaches a tooth length $2l$ to $x$. Do this for each connector in $\AA_2\setminus \AA_3$, and let $X'\subset X$ be the set of vertices attached to connectors in $A_3$. The vertices in $X'$ are attached to the connectors in $\AA_3$, and we will only find them teeth at the final stage of the embedding.

Now, take a subset $\BB_5\subset \BB_2\setminus \BB_4$ with $|\BB_5|=2|\AA_2|$. Using Lemma \ref{neat}, as before, we can find subsets $\BB_6\subset \BB_5$ and $\AA_4\subset \AA_3$, with $|\BB_6|\geq |\BB_5|-2m_2$ and $|\AA_4|\geq |\AA_3|-2m_2$, so that, for every $\UU\subset \BB_5\cup \AA_4$ with $|\UU|\leq 2m_2$, we have $|N_H(\UU,\BB_5\cup \AA_4)|\geq 4|\UU|$. Let $s_2:=\lfloor\mu_2 n/l\rfloor\leq s_1/12$ and pick a subset $\BB_7 \subset \BB_6$ so that $|\BB_7|=|\AA_3|-|\BB_4|-s_2$. This is possible, as $|\AA_3|\geq|\BB_4|+|\AA_3|/3\geq |\BB_4|+s_2$ and $|\BB_6|\geq 2|\AA_3|-m_2\geq |\AA_3|$.

Let $\CC=\AA_3$ and $\DD=\BB_4\cup \BB_7$. These sets of $(l,\gamma)$-connectors have the property that, given any subset $\UU\subset\CC$ with $|\UU|\leq m_2$, we have $|N_H(\UU^-,\DD^+)|\geq |N_H(\UU^-,\BB^+_4)|\geq 4|\UU|$ and, given any subset $\UU\subset \DD$ with $|\UU|\leq m_2$, if $\UU_1=\UU\cap \BB_4$ and $\UU_2=\UU\cap \BB_7$, then
\[
|N_H(\UU^+,\CC^-)|\geq \max_i|N_H(\UU^+_i,\CC^-)|\geq \max_i4|\UU_i|\geq 2|\UU|.
\]
We have embedded the tree $S$ apart from $|\CC|$ teeth length $2l$, where instead we have attached an $(l,\gamma)$-connector from $\CC$. We aim to embed the trees $T_1$ and $T_2$ without touching the connectors in $\CC\cup \DD$, before covering the remaining $s_2l$ vertices with $s_2$ paths of length $l-1$, whose endpoints have many neighbours among the sets $P^-$, $P\in \CC$. This will allow us to match these paths and the connectors indexed by $\DD$ onto the connectors indexed by $\CC$ and complete the embedding of $T$.
 
Reveal edges with probability $\e\log n/8n$ to get the graph $G_4$. Let $W$ be the set of vertices not in the partial embedding or in any of the connectors in $\CC$ or $\DD$. Let $L$ be the bipartite grouped graph with vertex sets $\CC^-$ and $W$ with respect to $G_4$. Edges are present independently with the same probability as in $K$, namely $p_2$. By Lemma \ref{mindegexp} and the choice of $\mu_2$, almost surely the graph $L$ has the $(4,D_2,\mu_2 n/l)$ property with $D_2=\e\mu_1\log\log n/ 800$. Call a vertex $v\in W$ \emph{good} if it has degree at least $D_2$ in the graph $L$, and \emph{bad} otherwise. Similarly to when we considered the graph $K$, as $D_2\geq p_3|\CC|/2$, by Proposition \ref{AKS1} there are almost surely $o(n)$ bad vertices.

Let $Z$ be the set of good vertices in $W$. The partial embedding of $S$ and the connectors cover $|S|-s_2l$ vertices, as they await $s_2$ paths of length $l-1$ to complete the embedding of $S$. As there are $o(n)$ bad vertices, $|Z|\geq n-|S|+s_2l-o(n)\geq (1-\e/14)n$. Reveal the final edges with probability $p_3=(1+\e/2)\log n/n$ to get the graph $G_5$. By Lemma \ref{HKS}, almost surely $\Delta(G_5)\leq 10\log n$  
and $d_{G_5}(u,Z)\geq \eta \log n$ for every vertex $u\in V(G)$, for $0<\eta =\eta(\e/14)<1/2$. By Proposition \ref{generalprops}, almost surely any two disjoint subsets of $V(G)$ of size $m_3=10n\log\log n / \log n$ have some edge between them in $G_5$. By Lemma \ref{mindegexp}, $G_5$ almost surely has the $(d_3,D_3,m_3)$-property for $d_3=\eta\mu_2\log n/ 10^5\log\log n$ and $D_3=\eta\mu_2\log n /3^6$.

Let $v_1$ and $v_2$ be the vertices in the embedding of $S'$ to which we embedded $t_1$ and $t_2$ respectively. By Lemma \ref{embedavoid}, we can in $G_5[W\cup\{v_1,v_2\}]$ embed $T_2$ and $T_3$ in $Z\cup\{v_1,v_2\}$ so that $t_1$ and $t_2$ are embedded on $v_1$ and $t_2$ respectively, and find $s_2$ disjoint paths $Q_i$ with length $l-1$ covering the vertices in $W$ not in these embeddings, so that the paths $Q_i$ have end vertices in $Z$. For each $i$, pick an end vertex of $Q_i$ and label it $q_i$, noting it is a good vertex.

We can now complete the embedding of $T$ by finding a matching in the bipartite grouped graph with vertex classes $\CC^-$ and $\DD^+\cup \mathcal{Q}$ with respect to the graph $G_2\cup G_4$, where $\mathcal{Q}=\{\{q_i\}:i\in[s_2]\}$. For a vertex $x\in X'$, let $P\in \CC$ be the connector with $P^+$ attached to $x$, and either let $q_j$ or $R$ be the vertex or connector matched to $P^-\in \CC^-$. Through $x$, $P$, and, respectively, either $Q_j$ or $R$, we can find a path length $2l$ to attach the required tooth to $x$. Doing this for each vertex $x\in X'$ completes the embedding of $T$. We will find such a matching by showing Hall's matching condition holds.

Let $\UU\subset \DD^+\cup \mathcal{Q}$ and suppose $|\UU|\leq 2m_2$. If $|\UU\cap \DD^+|\geq |\UU|/2$, then $|N_H(\UU\cap \DD^+,\CC)|\geq 2|\UU\cap \DD^+|\geq |\UU|$. If $|\UU\cap \mathcal{Q}|\geq |\UU|/2$, then, as $|\UU\cap \mathcal{Q}|\leq s_2$ and each vertex $q_i$ is good, we have, by the $(4,D_2,\mu_2n/ l)$-property of $L$, that $|N_L(\UU\cap \mathcal{Q})\geq 4|\UU\cap \mathcal{Q}|$. Hence, $|N_L(\UU\cap \mathcal{Q})|\geq |\UU|$, and Hall's condition is satisfied in this case.

If $|\CC|-m_2\geq |\UU|\geq 2m_2$, then $|\UU\cap \DD^+|\geq 2m_2-|\mathcal{Q}|\geq m_2$, so $|N_H(\UU\cap \DD^+,\CC)|\geq |\CC|-m_2\geq |\UU|$, and Hall's condition is also satisfied in this case.

Finally, if $|\UU|\geq |\CC|-m_2$, then, as above, we have $|N_H(\UU\cap \DD^+,\CC^-)|\geq |\CC|-m_2$. Hence, if $\VV=\CC^-\setminus N_H(\UU\cap \DD^+)$ then $|\VV|\leq m_2$. We have then $|N_H(\VV,\DD^+)|\geq |\VV|$. However, $N_H(\VV,\DD^+)\subset \DD^+\setminus \UU$, so $|\DD^+\setminus \UU|\leq|\VV|$. Therefore, $|N_H(\UU\cap \DD^+)|= |\DD|-|\VV|\geq|\DD|-|\DD^+\setminus \UU|\geq |\UU|$. Thus, Hall's condition is satisfied in all cases, and we can complete the embedding of $T$, as required.
\oof

\re When $10\leq k\leq 10^3\log n/\log\log n$, the only problem that arises in the previous proof is that the probability $p_2$ may not satisfy the upper bound in equation (\ref{bounds}). This is because if $l\leq 16\log n/\log\log n$, then the $(l,\gamma)$-connectors $P$ we use have $|P^+|=|P^-|=1$, where, for small $l$, this is much larger than $\gamma l$. Reducing the probability of an edge in $G_3$ and $G_4$ until the equation (\ref{bounds}) is satisfied fixes this problem.
\mark

In the above proof, if we increase $k$ then the probability $p_2$ will grow also, until it is sufficiently large that Lemma \ref{mindegexp} no longer holds. Therefore, we cannot find the matchings between vertices (in $X$ and $Q$) and the $(l,\gamma)$-connectors. This problem can be solved by replacing the vertices in $X$ and $Q$ with $(\log^6 n,\gamma)$-connectors to boost the edge probability for the matchings above the threshold for a matching to almost surely exist.

\pr[Sketch proof of Theorem \ref{combshort} when $k\geq \log^9 n$] Let $k_0=\lceil\log^6 n\rceil$, $l=\lfloor(k-k_0)/2\rfloor$, with $r=\lfloor l/3k_0\rfloor$, with $\mu_1$, $\mu_2$, $\gamma$ the same as before.

Using Lemma \ref{splittree}, divide the tree $T$ into the trees $S$, $T_1$ and $T_2$, with the same properties as before. Take $S$ and remove $s_1=\lfloor\mu_1 n/k\rfloor$ teeth length $l+(r+1)k_0$ to get $S'$.

Following the methods used in the proof of Theorem \ref{combshort} with $k_0$ in place of $k$, we can, by revealing edges with probability $\e\log n /8n$ three times, almost surely find a copy of $S'$ and $s_1$ $(k_0,\gamma)$-connectors, so that the following is true. Each vertex in $X$, the set of vertices which need teeth added, is attached to a different $(k_0,\gamma)$-connector $P_x$ using $P_x^+$. The only difference to the previous case is that $|X|$ is smaller, but this only makes it easier to find the matching and we can use Lemma \ref{mindegexp} as before.

In revealing these edges we can, as previously, almost surely get the expansion property we need to find $(l,\gamma)$-connectors using Lemma \ref{rotatingpaths}. Let $\AA$ be a set of $(3s_1-s_2)r$ $(k_0,\gamma)$-connectors, where $s_2=\lfloor\mu_2 n/k\rfloor$, and let $\BB$ be a set of $s_1$ $(l,\gamma)$ connectors, so that all the connectors are disjoint. Let $V$ be the vertices not in the copy of $S'$ or any of the connectors.

Reveal edges with the probability $\e\log n/ 32n$. Say that a vertex is \emph{good} if it has degree at least $\e\mu_1\log\log n/2400$ in the grouped graph $H$ with vertex classes $\AA^+$ and $V$, and \emph{bad} otherwise. The embedding of $S'$ and connectors we have found so far cover $|S|-s_2rk_0$ vertices. As previously, by revealing edges with probability $(1+\e/2)\log n /n$ and using Lemma \ref{embedavoid}, we may extend the embedding of $S'$ to cover $T_1$ and $T_2$, before splitting the remaining vertices in $V$ into $s_2r$ disjoint paths with length $k_0-1$ whose end vertices are good.

Taking the set of good end vertices of the paths we can match them into the connectors indexed by $A^+$ in $H$. This attaches a $(k_0,\gamma)$-connector to the end of each path, creating a $(3k_0,\gamma/3)$-connector.

We wish to turn the $3(s_1-s_2)r$ remaining $(k_0,\gamma)$-connectors in $\AA$ into $(3k_0,\gamma/3)$-connectors as well. This can be done by splitting them into three sets and revealing more edges between the sets with probability $\e\log n/32 n$. Almost surely, matchings will exist so that we may take groups of three connectors, one from each set, and join them end to end to create $(3k_0,\gamma/3)$-connectors.

We have embedded $S'$, $T_1$ and $T_2$, with a $(k_0,\gamma)$-connector attached to each of the $s_1$ vertices in need of a tooth. The remaining vertices are covered by $s_1$ $(l,\gamma)$-connectors and $rs_1$ $(3k_0,\gamma/3)$-connectors. We can finish the embedding very similarly to the embedding of $\mathrm{Comb}_{n,\sqrt{n}}$ in Section \ref{rotate}, working with the $(3k_0,\gamma/3)$-connectors as if they were vertices, but with an increased edge probability between them. Revealing more edges with probability $\e\log n/32n$, take a cycle through the $(3k_0,\gamma/3)$-connectors so that the cycle passes through each connector $P$, by coming in using $P^+$ and leaving using $P^-$. This is possible as the probability of an edge between two connectors is above the threshold for a directed Hamilton cycle in a random directed graph, as determined by McDiarmid~\cite{McD83}. Divide this cycle into $s_1$ sections of $r$ $(3k_0,\gamma/3)$-connectors. Revealing more edges with probability $\e\log n/32n$, find a matching to attach to each connector $P_x$, $x\in X$, some $(l,\gamma)$-connector, $Q_x\in\BB$ say, to which we connect a section of $r$ $(3k_0,\gamma/3)$-connectors. Taking a path through this structure allows us to attach a tooth the right length to each vertex $x\in X$, and complete the embedding of $T$.
\oof

\section{Proof of Theorem \ref{comblong}}\label{seclong}

Several sections of the proof of Theorem \ref{comblong} are very similar to those used in proving Theorem \ref{combshort}. Where there is a large overlap we will refer back to the proof of Theorem \ref{combshort} for the precise detail.

Theorem \ref{comblong} follows the outline of a theorem in \cite{selflove1}, where spanning trees were embedded in random graphs at a higher probability than the probability used here. Here we use connectors and work in the grouped graph, allowing the techniques to be applied with a lower probability.

\pr[Proof of Theorem \ref{comblong}] We will prove the theorem with $\a=\e/7$; the general case follows by reducing either $\a$ or $\e$ until this holds. We will reveal edges in stages to get the graphs $G_1$, $G_2$ and $G_3$. Because Theorem \ref{comblong} is a universal result, we wait until the final graph $G=G_1\cup G_2\cup G_3$ is fully revealed before embedding an arbitrary tree with many bare paths, using the properties of $G$. We will reveal edges independently with probability $\e \log n/4n$ in $G_1$ and $G_2$, and with probability $(1+\e/2)\log n/ n$ in $G_3$, so that we may compare $G$ with $\GG(n,(1+\e)\log n/ n)$.

Let $k=\log^9 n$, $l_1=\lceil\log^6 n/12\rceil$, $l_2=\lceil\log^3 n\rceil$ and $l_0=6l_2l_1+2$, so that $k\geq l_0\geq k/2$. Let $\mu$ be sufficiently small that, by Lemma \ref{mindegexp}, a random graph $\GG(n,p)$, for any $p=p(n)\leq \log^{10}n/ n$, almost surely has the $(4,\e^2\log\log n/10^5,\mu n\log\log n/p)$-property. Let $s_1=\lfloor \e n/14l_0\rfloor$ and $s_2=\lfloor\mu n/l_1\rfloor$, so that $s_2<s_1l_2/2$.

Revealing the edges of $G_1$, we can almost surely find the following sets of disjoint connectors in $G_1$, using Lemma \ref{rotatingpaths}, as previously. Letting $\gamma=\log\log n/ 16\log n$, the sets $\AA$ and $\BB$, contain respectively $(s_1l_2-s_2)$ and $(2s_1l_2+2s_1)$ $(l_1,\gamma)$-connectors, and the set $\CC$ contains $s_1l_2$ $(3l_1,\gamma)$-connectors. In total, these connectors cover $s_1l_0-s_2l_1\leq \e/14 -\mu n/2$ vertices.

Reveal the edges of $G_2$. Let $H$ be the bipartite grouped graph with vertex classes $\BB^+$ and $\AA^+\cup \AA^-$ with respect to $G_2$. Then $H$ has edges present independently with probability $p_1$, where
\[
p_1\geq \frac{\e(\gamma l_1)^2\log n}{8n}\geq \frac{\e(\log n)^{11}(\log\log n)^2}{2^{20}n}\geq \frac{10^5\e\log^5n \log\log n}{|\AA|}.
\]
Almost surely, by Proposition \ref{generalprops}, there will be an edge in $H$ between any two disjoint vertex sets if they each contain at least $m_1=\lceil|A|/10\log^5n\rceil$ vertices. Moreover, almost surely, for every subset $U\subset V(H)$, with $|U|\leq |\AA|/10\log^5 n$, we have $|N_K(U)|\geq \log^5 n|U|$, by Lemma \ref{mindegcond}.

Let $K$ be the grouped graph with vertex set $\BB^-\cup \CC^+\cup \CC^-$ with respect to $G_2$. The probability an edge in $K$ is present is at least $p_1$. As for the graph $H$, almost surely there will be an edge in $K$ between any two disjoint sets containing at least $m_1$ vertices each. By Lemma \ref{mindegcond}, as $m_1\leq |\CC|/2\log^5n$, almost surely any subset $U\subset V(K)$ with $|U|\leq m_1$ must have $|N_K(\UU,\CC^+)|,|N_K(\UU,\CC^-)|\geq |\UU|\log^5 n$.

Let $V$ be the vertices of $G$ not in any of the connectors. Let $L$ be the bipartite grouped graph with vertex sets $B^+$ and $V$ with respect to $G_2$. Then $L$ has edges present independently with some probability $p_2$, where
\[
\frac{\e l_1\log\log n}{128n}\leq p_2=1-\left(1-\frac{\e\log n}{8n}\right)^{\left\lceil\frac{l_1\log\log n}{16\log n}\right\rceil}\leq \frac{\e l_1\log\log n}{n}.
\]
Let $D_1=\e^2\log\log n/10^5$. Almost surely, the graph $L$ will the $(4,D_1,\mu n/l_1)$-property by Lemma \ref{mindegexp} and the choice of $\mu$.

Let a vertex $v\in V$ be \emph{good} if it has at least $D_1$ neighbours in $L$, and \emph{bad} otherwise. Let $V_1$ contain these good vertices. As previously, because $D_1\leq p_3|B|/2$, by Proposition \ref{AKS1} we almost surely have $o(n)$ bad vertices. Therefore $|V_1|=|V|-o(n)\geq (1-\e/14)n+\mu n/2-o(1)\geq (1-\e/14)n$. Reveal edges with probability $(1+\e/2)\log n/ n$ to get the graph $G_3$. By Lemma \ref{HKS}, almost surely $\Delta(G_3)\leq 10\log n$  
and $d_{G_3}(u,V_1)\geq \eta \log n$ for every $u\in V(G)$, for $0<\eta =\eta(\e/14)<1/2$. By Proposition \ref{generalprops}, almost surely any two sets of size at least $m_2=\lceil 10\log\log n/p_2\rceil$ have some edge between them in $G_3$. By Lemma \ref{mindegexp}, $G_3$ will almost surely have the $(d_2,D_2,m_2)$-property for $d_2=\eta\mu \log n/ 3^{11}\log\log n$ and $D_2=\eta\mu\log n/ 3^6$. 

The graph $G=G_1\cup G_2\cup G_3$ is now fixed. Let $T$ be any tree in $\TT(n,\Delta)$ with at least $\e n/14 k$ disjoint bare paths length $k$. Remove $s_1$ bare paths length $l_0+1$ from the tree $T$ to get a forest $T'$ with $n-s_1l_0$ vertices. Due to the vertices covered by connectors, $|V|=n-s_1l_0+s_2l_1$. By adding dummy edges we can turn the forest $T'$ into two trees, $T_1$ and $T_2$ say. Using Lemma \ref{embedavoid}, as $|V_1|=|V(T')|+s_2l_1-o(n)$, we may embed the trees $T_1$ and $T_2$ into the good vertices of the graph and divide the rest of the vertices in $V$ into $s_2$ paths length $l_1-1$, each with two good ends. Let these paths be $R_i$, $i\in[s_2]$, each having end vertices $r^+_i$ and $r^-_i$.

Suppose that to make the embedding of $T'$ into an embedding of $T$ we have to connect, for each $i\in[s_1]$, $x_i$ to $y_i$ with a bare path of length $l_0+1$. Let $R^+=\{r_i^+:i\in[s_2]\}$, $R^-=\{r_i^-:i\in[s_2]\}$, $X=\{x_i:i\in[s_1]\}$ and $Y=\{y_i:i\in[s_1]\}$.

We will find a matching between the set $\BB^+$ and the set $F=\AA^+\cup \AA^-\cup R^+\cup R^-\cup X\cup Y$, so that if two sets are matched together then there is an edge between them in either $H$ or $L$. Note that $|F|=2(s_1l_2-s_2)+2s_1+2s_2=|\BB^+|$. Letting $\mathcal{Q}=R^+\cup R^-\cup X\cup Y$, $\DD=\AA^+\cup\AA^-$, and $\CC=\BB^+$, we have an extremely similar situation to the matching found at the end of the proof of Theorem \ref{combshort}, and by similarly showing Hall's matching condition holds, we can demonstrate that such a matching exists.

This matching allows us to attach an $(l_1,\gamma)$-connector to each vertex in $X$ and $Y$, attach an $(l_1,\gamma)$-connector to each end of the path $R_i$ so that it becomes a $(3l_1,\gamma/3)$-connector, and attach an $(l_1,\gamma)$-connector to each end of the $(l_1,\gamma)$-connectors in $\AA$ so that they too become $(3l_1,\gamma/3)$-connectors.

In summary, we have an embedding of $T'$ with an $(l_1,\gamma)$-connector connected to each vertex in $X$ and $Y$. The remaining vertices have been divided into $s_1l_2$ $(3l_1,\gamma)$-connectors in $\CC$ and $s_1l_2$ $(3l_1,\gamma/3)$-connectors formed from the matching.

Let $\mathcal{E}$ be the set of $(3l_1,\gamma/3)$ connectors we have formed, and let $x_i$ be matched to $X_i\in \BB$ and $y_i$ be matched to $Y_i\in \BB$, for each $i\in[s_1]$. Form a new directed graph, $M$, on the vertex set $\mathcal{X}=\CC\cup\mathcal{E}\cup\{X_i,Y_i:i\in[s_1]\}$. For each pair of connectors $P,Q\in \mathcal{X}$, if there is an edge from $P^+$ to $Q^-$ in the graph $G$, then let $\vec{PQ}$ be an edge of $M$. Let $\mathcal{W}=\CC\cup\mathcal{E}$. By the properties of the graph $K$, we can see that the graph $M$ with pairs $(X_i,Y_i)$, $i\in[s_1]$, satisfies the properties of Theorem \ref{pathcoverexpander} with the set $\mathcal{W}$. Therefore, we can cover $M$ with directed paths of length $2s_2$ which are disjoint and connect the pairs $(X_i,Y_i)$.

For each $i\in[s_1]$, start from $x_i$ and find a path length $l_0+1$ to $y_i$ which passes through the connectors in the order determined by the directed path in $M$ from $X_i$ to $Y_i$.  Doing this creates disjoint paths of length $l_0+1$ between all the pairs of vertices $(x_i,y_i)$, as required, and thus completes the embedding of $T$.
\oof

\begin{acknowledgements} The author would like to thank Andrew Thomason for his help and suggestions.
\end{acknowledgements}

\bibliographystyle{plain}
\bibliography{rhmreferences}

\begin{thebibliography}{10}

\bibitem{AKS07}
N.~Alon, M.~Krivelevich, and B.~Sudakov.
\newblock Embedding nearly-spanning bounded degree trees.
\newblock {\em Combinatorica}, 27(6):629--644, 2007.

\bibitem{BCPS10}
J.~Balogh, B.~Csaba, M.~Pei, and W.~Samotij.
\newblock Large bounded degree trees in expanding graphs.
\newblock {\em Electronic Journal of Combinatorics}, 17(1):R6, 2010.

\bibitem{bollo1}
B.~Bollob{\'a}s.
\newblock {\em Random graphs}.
\newblock Springer, 1998.

\bibitem{PH01}
P.E. Haxell.
\newblock Tree embeddings.
\newblock {\em Journal of Graph Theory}, 36(3):121--130, 2001.

\bibitem{HKS09}
D.~Hefetz, M.~Krivelevich, and T.~Szab{\'o}.
\newblock Hamilton cycles in highly connected and expanding graphs.
\newblock {\em Combinatorica}, 29(5):547--568, 2009.

\bibitem{HKS12}
D.~Hefetz, M.~Krivelevich, and T.~Szab{\'o}.
\newblock Sharp threshold for the appearance of certain spanning trees in
  random graphs.
\newblock {\em Random Structures and Algorithms}, 41(4):391--412, 2012.

\bibitem{JKS12}
D.~Johannsen, M.~Krivelevich, and W.~Samotij.
\newblock Expanders are universal for the class of all spanning trees.
\newblock In {\em Proceedings of the Twenty-Third Annual ACM-SIAM Symposium on
  Discrete Algorithms}, pages 1539--1551. SIAM, 2012.

\bibitem{KLW14a}
J.~Kahn, E.~Lubetzky, and N.~Wormald.
\newblock Cycle factors and renewal theory.
\newblock {\em arXiv preprint arXiv:1401.2707}, 2014.

\bibitem{KLW14b}
J.~Kahn, E.~Lubetzky, and N.~Wormald.
\newblock The threshold for combs in random graphs.
\newblock {\em arXiv preprint arXiv:1401.2710}, 2014.

\bibitem{MK10}
M.~Krivelevich.
\newblock Embedding spanning trees in random graphs.
\newblock {\em SIAM Journal on Discrete Mathematics}, 24(4):1495--1500, 2010.

\bibitem{McD83}
C.~McDiarmid.
\newblock General first-passage percolation.
\newblock {\em Advances in Applied Probability}, 15(1):pp. 149--161, 1983.

\bibitem{selflove1}
R.H. Montgomery.
\newblock Embedding bounded degree spanning trees in random graphs.
\newblock {\em Preprint}, 2014.

\bibitem{posa76}
L~P{\'o}sa.
\newblock Hamiltonian circuits in random graphs.
\newblock {\em Discrete Mathematics}, 14(4):359--364, 1976.

\end{thebibliography}

\end{document}